\documentclass[10pt,a4paper]{amsart}

 \usepackage{amssymb,amsmath,amsfonts,amsthm} 

\theoremstyle{plain}
\newtheorem{thm}{Theorem}[section]
\newtheorem{prop}{Proposition}[section]
\newtheorem{lem}{Lemma}[section]
\newtheorem{cor}{Corollary}[section]

\theoremstyle{definition}
\newtheorem{defi}{Definition}[section]
\newtheorem*{nota}{Notation}

\theoremstyle{remark}
\newtheorem*{rem}{Remark}
\newcommand{\G}{\mathcal{G}}

\begin{document}
\title{On deformations of the filiform Lie superalgebra $L_{n,m}$}
\author{M. GILG}
\address{Marc GILG, Universit\'e de Haute-Alsace,
Laboratoire de Math\'ematiques,
4 rue des Fr\`eres Lumi\`ere,
68 093 Mulhouse Cedex, France}
\email{M.Gilg@univ-mulhouse.fr}

\date{}
\maketitle

\section*{Abstract}

Many work was done for filiform Lie algebras
defined by M. Vergne \cite{ver}. An interesting fact
is that this algebras are obtained by deformations of
the filiform Lie algebra $L_{n,m}$. This was used for
classifications in \cite{kj}.
Like filiform Lie algebras, filiform Lie
superalgebras are obtained by nilpotent
deformations of the Lie superalgebra $L_{n,m}$.
In this paper, we recall this fact and we 
study even cocycles of the superalgebra $L_{n,m}$
which give this nilpotent deformations. 
A family of independent bilinear maps will help us to
describe this cocycles. At the end
an evaluation of the dimension of the space $Z_0^2(L_{n,m},L_{n,m})$
is established. The description of this cocycles can help us to get
some classifications which was done in \cite{mg,mg1}. 

\section{Deformation of Lie superalgebras}

\subsection{Nilpotent Lie superalgebras}

\begin{defi} A  $\mathbb{Z}_2$-graded vector space $\G=\G_0\oplus \G_1$ over an algebraic
closed field is a
Lie superalgebra if there exists a bilinear product $[,]$ over $\G$ such
that 
$$\begin{aligned}
\ &[G_{\alpha},\G_{\beta}]\subset \G_{\alpha+\beta\mod 2},\\ 
&[g_{\alpha}, g_{\beta}]=(-1)^{\alpha.\beta} [g_{\beta},g_{\alpha}]
\end{aligned}$$
 for all
$g_{\alpha}\in\G_{\alpha}$ and $g_{\beta}\in\G_{\beta}$ and satisfying 
Jacobi identity:
$$(-1)^{\gamma .\alpha}[A,[B,C]]+(-1)^{\alpha .\beta}[B,[C,A]]+(-1)^{\beta .\gamma}[C,[A,B]]=0$$
for all $A\in\mathcal{G}_\alpha$, $B\in\mathcal{G}_\beta$ and $C\in\mathcal{G}_\gamma$.
\end{defi}

For such a Lie superalgebra we define the
lower central series
$$\left\{
\begin{aligned}
&C^0(\mathcal{G})=\mathcal{G},\\
&C^{i+1}(\mathcal{G})=[\mathcal{G},C^i(\mathcal{G})].
\end{aligned}
\right.$$

\begin{defi}
A Lie superalgebra $\mathcal{G}$ is nilpotent if there exist
an integer $n$ such that $C^n(\mathcal{G})=\{0\}$.
\end{defi}

We define for a Lie superalgebra 
\label{seq}
$\mathcal{G}=\mathcal{G}_0\oplus\mathcal{G}_1$ two sequences :
$$C^0(\mathcal{G}_0)=\mathcal{G}_0, \qquad C^{i+1}(\mathcal{G}_0)=[\mathcal{G}_0,C^i(\mathcal{G}_0)]$$
and $$C^0(\mathcal{G}_1)=\mathcal{G}_1,\qquad C^{i+1}(\mathcal{G}_1)=[\mathcal{G}_0,C^i(\mathcal{G}_1)]$$

\begin{thm}
Let $\mathcal{G}=\mathcal{G}_0\oplus\mathcal{G}_1$ be a Lie superalgebras.
Then $\mathcal{G}$ is nilpotent if and only if there exist $(p,q)\in\mathbb{N}^2$ such that
$C^p(\mathcal{G}_0)=\{0\}$ and $C^q(\mathcal{G}_1)=\{0\}$.
\end{thm}

\begin{proof}
If the Lie superalgebra $\G=\G_0\oplus \G_1$ is nilpotent
the existence of $(p,q)$ such that
$C^p(\mathcal{G}_0)=\{0\}$ and $C^q(\mathcal{G}_1)=\{0\}$
is obvious.

For the converse, assume that there exist $(p,q)$ such that
$C^p(\mathcal{G}_0)=\{0\}$ and $C^q(\mathcal{G}_1)=\{0\}$,
then every operator  $ad(X)$ with $X\in\G_0$ is
nilpotent.  
Let  $Y\in \mathcal{G}_1$, as 
$$ad(Y)\circ ad(Y)=\frac{1}{2} ad([Y,Y])$$
$[Y,Y]$ is an element of $\G_0$, then $ad([Y,Y])$ is nilpotent.
This implies that $ad(Y)$ is nilpotent for every $Y\in \G_1$.
By Engel's theorem for Lie superalgebras \cite{sch}, this implies that $\mathcal{G}$
is nilpotent Lie superalgebra.
\end{proof}

\begin{defi}
Let $\G$ be a nilpotent Lie superalgebra, the super-nilindex of $\G$ is the 
pair $(p,q)$ such that :
$C^p(\mathcal{G}_0)=\{0\}$, $C^{p-1}(\mathcal{G}_0)\ne\{0\}$
and $C^q(\mathcal{G}_1)=\{0\}$, $C^{q-1}(\mathcal{G}_1)\ne\{0\}$. It is and invariant
up to isomorphism.
\end{defi}

\subsection{Cohomology}

We recall some definition from \cite{fu}.

By definition, the superspace of $q$-dimensional cocycles of the
Lie superalgebra $\G=\G_0\oplus\G_1$ with coefficient in the $\G$-module
$A=A_0\oplus A_1$ is given by

$$C^q(\G ;A)=\bigoplus_{q_0+q_1=q} Hom( \bigwedge^{q_0}\G_0\otimes\bigvee^{q_1}\G_1,A).$$
This space is graded by $C^q(\G ;A)=C_0^q(\G ;A)\oplus C_0^q(\G ;A)$ with
\[
C_p^q(\G ;A)=\bigoplus_{\substack{q_0+q_1=q\\ q_1+r\equiv p\mod 2}} 
Hom( \bigwedge^{q_0}\G_0\otimes\bigvee^{q_1}\G_1,A_r)
\]

The differential 
$$d:\quad C^q(\G;A)\longrightarrow C^{q+1}(\G;A)$$
defined by
$$\begin{aligned}
&\text{$d\ c(u_1,\dots ,u_{q_0},v_1,\dots ,v_{q_1})$}\\ 
&=\text{$\sum_{1\leq s<t\leq q_0} (-1)^{s+t-1}\ c([u_s,u_t],u_1,\dots ,\hat{g}_s,\dots ,\hat{g}_t,\dots ,u_{q_0},v_1,\dots,v_{q_1})$}\\
&\quad+\text{$\sum_{s=1}^{q_0}\sum_{t=1}^{q_1} (-1)^{s-1}\ c(u_1,\dots,\hat{g}_s,\dots u_{q_0},[u_s,v_t],v_1,\dots ,\hat{h}_t,\dots ,v_{q_1})$}\\
&\quad+ \text{$\sum_{1\leq s<t\leq q_1} c([v_s,v_t],u_1,\dots ,u_{q_0},v_1,\dots,\hat{h}_s,\dots ,\hat{h}_t,\dots v_{q_1})$}\\
&\quad+\text{$\sum_{s=1}^{q_0} (-1)^s\ u_s\ c(u_1,\dots ,\hat{g}_s,\dots u_{q_0},v_1,\dots v_{q_1})$}\\
&\quad+(-1)^{q_0-1}\text{$\sum_{s=1}^{q_1} v_s\ c(u_1,\dots u_{q_0},v_1,\dots ,\hat{h}_s,\dots, v_{q_1})$}\\
\end{aligned}$$
where $c\in C^q(\G;A)$, $u_1,\dots ,u_{q_0}\in \G_0$ and $v_1,\dots ,v_{q_1}\in \G_1$ and satisfies 
$$\left\{
\begin{aligned}
&d\circ d=0\\
&d(C_p^q(\G;A))\subset C_p^{q+1}(\G;A).
\end{aligned}
\right.$$
for $q=0,1,2\dots$ and $p=0,1$.

Let be  $d_j:\ C_p^{j}(\G;A)\longrightarrow C_p^{j+1}(\G;A)$ with $p=0$ or $p=1$ the restriction
of $d$ to the space $C_p^{j}(\G;A)$. This operator permit to define the spaces :
$$H_p^j(\G,A)=\frac{Z_p^j(\G;A)}{B_p^j(\G;A)}$$
where $p=0$ or $p=1$.
Therefore we have :
\begin{itemize}
\item $Z^j(\G,A)=Z_0^j(\G,A)\oplus Z_1^j(\G,A)$. 
\item $B^j(\G,A)=H_0^j(\G,A)\oplus B_1^j(\G,A)$. 
\item $H^j(\G,A)=H_0^j(\G,A)\oplus H_1^j(\G,A)$.    
\end{itemize}

\subsection{Algebraic variety of nilpotent Lie superalgebras}

We recall some facts from \cite{goz}.

Let $\mathcal{L}^{n}_{p,q}$ be the set of Lie superalgebras law
over $\mathbb{C}^{n}=\mathbb{C}^{p+1}\oplus \mathbb{C}^{q}$. 
Let $(X_{1},X_{2},\ldots X_{p+1},Y_{1},Y_{2},\ldots Y_{q})$ be a graded base of it.
For $\mu \in \mathcal{L}^{n}_{p,q}$ we set  :
$$\left\{
\begin{aligned}
\mu (X_{i},X_{j})&=\sum _{k=1}^{p+1}C^{k}_{i,j}\ X_{k}\qquad 1\leq i<j\leq p\\
\mu (X_{i},Y_{j})&=\sum _{k=1}^{q}D^{k}_{i,j}\ Y_{k}\qquad 1\leq i\leq p+1,\; 1\leq j\leq q\\
\mu (Y_{i},Y_{j})&=\sum _{k=1}^{p+1}E^{k}_{i,j}\ X_{k\qquad 1\leq i\leq j\leq q}\\
\end{aligned}
\right.$$
with $C^{k}_{j,i}=-C^{k}_{i,j}$ and $E^{k}_{j,i}=E^{k}_{i,j}$.

The elements $\{ C_{i,j}^k,D_{i,j}^k,E_{i,j}^k\}_{i,j,k}$ are called
{\it structure constants} of the Lie superalgebra with respect to basis
$(X_1,\dots ,X_{p+1},Y_1,\dots ,Y_q)$.
\label{salp}
The Jacobi identities 
show that  $\mathcal{L}^{n}_{p,q}$ is an algebraic sub-variety of
$\mathbb{C}^N$ with
 $$N=(p+1)^{2}(\frac{p}{2})+2(p+1)q^{2}.$$

Let $V=V_0\oplus V_1$ be a $\mathbb{Z}_2$-graded vector space
of dimension $n$ with $\dim V_0=p+1$ and $\dim V_1=q$.
Let $G(V)$ be the group of linear map  of the type $g=g_0+g_1$
where $g_0\in GL(V_0)$ and $g_1\in GL(V_1)$. This group is
isomorphic to $GL(V_0)\times GL(V_1)$. 

The algebraic group $G(V)$ acts on the variety 
$\mathcal{L}^n_{p,q}$ in the following way  :
\[
(g.\phi)(x,y)=g_{\alpha+\beta}
(\phi(g_{\alpha}^{-1}(a),g_{\beta}^{-1}(b)))\quad \forall a\in V_{\alpha},\ \forall b\in V_{\beta},
\]
with $g\in G(V)$ and $\phi\in \mathcal{L}_{p,q}^n$.

\subsection{Deformations of Lie superalgebras}

Let $\mathcal{G}$ be a Lie superalgebra over a field $k$, $V$ be the underlying
vector space and $\nu_0$ the law of $\mathcal{G}$.
Let $k\left[ \left[t\right] \right] $ be the power series ring in one variable $t$. 
Let $V\left[ \left[ t\right] \right] $ be the $k\left[ \left[ t\right] \right]$-module
$V\left[ \left[ t\right] \right] =V\otimes _kk\left[ \left[t\right] \right] $.  
One can obtain an extension of $V$ with a structure of vector
space by extending the coefficient domain from $k$ to $k\left( \left(
t\right) \right) $, the quotient power series field of $k\left[ \left[
t\right] \right] $.Any bilinear map $f:V\times V\rightarrow V$ (in
particular the multiplication in $\G$) can be extended to a bilinear map from 
$V\left[ \left[ t\right] \right] \times V\left[ \left[ t\right] \right] $ to 
$V\left[ \left[ t\right] \right] .$

\begin{nota}
Let $A_{p,q}^2$ be the set of bilinear forms  
$$\phi:\ k^n\times k^n\rightarrow k^n$$
satisfying :
$$\left\{
\begin{aligned}
&\phi(V_i,V_j)\subset V_{i+j\mod 2},\\
&\phi(v_i,v_j)=(-1)^{d(v_i).d(v_j)}\ \phi(v_j,v_i).
\end{aligned}
\right.$$

where $k^n=V_0\oplus V_1$ is a  $\mathbb{Z}_2$-graded vector space ,
$\dim V_0=p+1$, $\dim V_1=q$ and $v_i\in V_{d(v_i)}$, $v_i\in V_{d(v_i)}$.
\end{nota}

\begin{defi} Let $\nu_0$ be the law of the
Lie superalgebra $\G$. A deformation of $\nu_0$ is a one
parameter family $\nu_t$ in $k\left[ \left[ t\right] \right] \otimes V$.

$$\nu_t=\nu_0+t.\nu_1+t^2.\nu_2+\dots$$
where $\nu_i\in A_{p,q}^2$ for $i\geq 1$,
$\nu_t$ satisfy the Jacobi formal identities : 
$$(-1)^{\gamma .\alpha}\nu_t(A,\nu_t(B,C))+(-1)^{\alpha .\beta}\nu_t(B,\nu_t(C,A))+(-1)^{\beta .\gamma}\nu_t(C,\nu_t(A,B))=0,$$
For all $A\in\mathcal{G}_\alpha$, $B\in\mathcal{G}_\beta$ and $C\in\mathcal{G}_\gamma$.
\end{defi}

The coefficient of $t^k$ of the formal Jacobi identity is
$$\left\{ 
\begin{aligned}
&\sum_{i=0}^k (-1)^{\gamma .\alpha}\nu_i(A,\nu_{k-i}(B,C))+
(-1)^{\alpha .\beta}\nu_i(B,\nu_{k-i}(C,A))+(-1)^{\beta .\gamma}\nu_i(C,\nu_{k-i}(A,B))=0\\
&\quad k=0,1,2,\text{$\dots$}
\end{aligned}
\right.$$
for all $A\in\mathcal{G}_\alpha$, $B\in\mathcal{G}_\beta$ and $C\in\mathcal{G}_\gamma$.
This last relations are called the deformation equations.

For $k=0$ we get the Jacobi identity of the Lie superalgebra  $\nu_0$. 
For $k=1$ the condition on the coefficient $t$ implies the next proposition  :
\begin{prop}
\label{def}
Let $\nu_0$ be a Lie superalgebra and 
$\nu_t$ of it  :
$$\nu_t=\nu_0+t.\nu_1+t^2.\nu_2+\dots$$
then $\nu_1$ is an even 2-cocycle of the
Lie superalgebra $\nu_0$ ($\nu_1\in Z_0^2(\nu_0,\nu_0)$).
\end{prop}

\subsection{Deformation in $\mathcal{N}_{n,m}^{p,q}$ }

Let $\G=\G_0\oplus\G_1$ be a nilpotent Lie superalgebra of $\mathcal{N}_{n,m}^{p,q}$ with 
multiplication  $\nu_0$ and $\nu_t$ be a deformation of it.\\
We write  $\nu_0=\mu_0+\rho_0+b_0$ where :
$$\begin{aligned}
\mu_0&\in \hom (\mathcal{G}_0\wedge \mathcal{G}_0,\mathcal{G}_0)\\
\rho_0&\in \hom(\mathcal{G}_0\otimes\mathcal{G}_1,\mathcal{G}_1)\\
b_0&\in \hom (\mathcal{G}_1\vee\mathcal{G}_1,\mathcal{G}_0)
\end{aligned}$$
For $\nu_t$  to be a deformation in $\mathcal{N}_{n,m}^{p,q}$, we must have :
$$\text{(N)}\quad
\left\{
\begin{aligned}
&\text{$\nu_t(x_1,\nu_t(x_1,\dots ,\nu_t(x_p,x_0)\dots))=0$}\\
&\text{$\nu_t(x_1,\nu_t(x_1,\dots ,\nu_t(x_q,y)\dots))=0$}\\
\end{aligned}
\right.$$
for all $x_i$ in $\G_0$  and  $y$ in $\G_1$.
Proposition \ref{def} implies that  $\nu_1\in Z_0^2(\nu_0,\nu_0)$. 
Let be $\nu_1=\psi_1+\rho_1+b_1$ with :
$$\begin{aligned}
\psi_1&\in \hom (\mathcal{G}_0\wedge \mathcal{G}_0,\mathcal{G}_0)\\
\rho_1&\in \hom(\mathcal{G}_0\otimes\mathcal{G}_1,\mathcal{G}_1)\\
b_1&\in \hom (\mathcal{G}_1\vee\mathcal{G}_1,\mathcal{G}_0)
\end{aligned}$$

\section{Filiform Lie superalgebras}

\begin{defi}
Let $\mathcal{G}=\mathcal{G}_0\oplus\mathcal{G}_1$ be a nilpotent Lie superalgebra with
$\dim\mathcal{G}_0=n+1$ and $\dim\mathcal{G}_1=m$. $\G$ is called filiform if
it's super-nilindex is $(n,m)$.
 
We will note $\mathcal{F}_{n,m}$ the set of
filiform Lie superalgebras.
\end{defi}

\begin{rem}
We can write 
the set of filiform Lie superalgebras as the
complement of the closed set  for the Zariski topology of the nilpotent
superalgebras with super-nilindex $(k,p)$ such that $k\leq n-1$ and $p\leq m-1$.
Hence the set of filiform Lie superalgebras is an open set of the variety of
nilpotent Lie superalgebras.
\end{rem}

As for the filiform Lie algebras \cite{ver}, there exists an adapted base of a filiform Lie
superalgebra :

\begin{thm}
\label{adb}
Let $\mathcal{G}=\mathcal{G}_0\oplus\mathcal{G}_1$ be a filiform Lie superalgebra with
$\dim\mathcal{G}_0=n+1$ and $\dim\mathcal{G}_1=m$.
Then there exists a base $\{X_0,X_1,\dots X_n,Y_1,Y_2,\dots Y_m\}$ of $\G$
with $X_i\in \G_0$ and $Y_i\in\mathcal{G}_1$
such that :
$$\left\{
\begin{aligned}
\ &[X_0,X_i]=X_{i+1}\quad 1\leq i\leq n-1, \quad [X_0,X_n]=0;\\
&[X_1,X_2]\in \mathbb{C}.X_4+\mathbb{C}.X_5+\dots + \mathbb{C}.X_n;\\
&[X_0,Y_i]=Y_{i+1}\quad 1\leq i\leq m-1,\quad [X_0,Y_m]=0.
\end{aligned}
\right.$$
\end{thm}
The proof is the same as for Lie algebras \cite{ver}( see also \cite{mg1}).

{\bf Example} Define the superalgebra $L_{n,m}=L_{n,m}^0\oplus L_{n,m}^1$ by
$$\left\{
\begin{aligned}
\ &[X_0,X_i]=X_{i+1}\quad 1\leq i\leq n-1,\\
&[X_0,Y_i]=Y_{i+1}\quad 1\leq i\leq m-1.\\
\end{aligned}
\right.$$
where the other brackets vanished,
$\dim L_{n,m}^0=n+1$, $\dim L_{n,m}^1=m$ and\\
$\{ X_0,X_1,\dots X_n,Y_1,\dots, Y_m\}$ is an adapted base.
The law of $L_{n,m}$ is written by $\mu=\mu_0+\rho_0$ where
$\mu_0$ is the law of the Lie algebras $L_{n,m}^0$ and $\rho_0$ is
the representation associated to the $L_{n,m}^0$-module $L_{n,m}^1$.

\begin{prop}
\label{pdef}
Every filiform Lie superalgebra $\G=\G_0\oplus\G_1$ such that
$\dim\G_0=n+1$ and $\dim\G_1=m$ can be written :
$$[\bullet,\bullet ]=\mu_0+\rho_0+\Phi$$
with $\Phi$ satisfying :

$$\begin{aligned}
&\Phi\in Z_0^2(L_{n,m},L_{n,m})\\
&\Phi(X_0,Z)=0 \quad \forall Z\in\G\\
&\Phi(S_i^0,S_j^0)\subset S_{i+j+1}^0\text{ if }i+j<n\\
&\Phi(X_i,X_{n-i})=(-1)^i\alpha\ X_n\text{ where }
\alpha=0\text{ if $n$ is even }\\
&\Phi(S_i^0,S_j^1)\subset S_{i+j}^1
\end{aligned}$$
where $\mu_0+\rho_0$ is the law of $L_{n,m}$, and $S_.^0$, $S_.^1$
the filtrations associated to the graduations of $\mathcal{G}_0$
and $\mathcal{G}_1$ by the sequences given in \ref{seq}.

\end{prop}

\begin{proof}

Using the theorem \ref{adb}, for every filiform lie superalgebra we have
an adapted base $\mathcal{B}$ $\{ X_0,X_1,\dots , X_n,Y_1,Y_2,\dots Y_n\}$
 such that the product of $\G$ is given by :
$$[\bullet,\bullet ]=\mu_0+\rho_0+\Phi$$
where $\Phi[X_0,Z]=0$ for every vector $Z\in\G$.
This product satisfies the Jacobi identity 
\begin{equation}
\label{jacop1}
(\mu_0+\rho_0)\circ \Phi+\Phi\circ(\mu_0+\rho_0)+\Phi\circ\Phi=0.
\end{equation}
Let $Z_i,Z_j$ be two vectors of the adapted base $\mathcal{B}$ of $\G$.
We have $\Phi\circ\Phi(X_0,Z_i,Z_j)=0$ because $\Phi(X_0,\bullet)=0$.
The relation (\ref{jacop1}) becomes :
$$((\mu_0+\rho_0)\circ \Phi+\Phi\circ(\mu_0+\rho_0))(X_0,Z_i,Z_j)=0$$
Also  we have for every $Z_i,Z_j,Z_k\in\mathcal{B}\setminus \{ X_0\}$ :
$$((\mu_0+\rho_0)\circ \Phi+\Phi\circ(\mu_0+\rho_0))(Z_i,Z_j,Z_k)=0.$$
As the superalgebra is nilpotent, $X_0\notin Im \Phi$ shows : 
$$((\mu_0+\rho_0)\circ \Phi+\Phi\circ(\mu_0+\rho_0))(U,V,W)=0.$$
for all $U,V,W\in\mathcal{G}$, where $\circ$ is the graded
Nijenhuis-Richardson bracket.
This implies that $\Phi\in Z_0^2(L_{n,m},L_{n,m})$.
The filtrations :
$$\left\{
\begin{aligned}
\ [S_i^0,S_j^0]\subset S_{i+j}^0\\
[S_i^0,S_j^1]\subset S_{i+j}^1\\
\end{aligned}
\right.$$
associated to the graduations $C^{i}(\mathcal{G}_0)$ and $C^{i}(\mathcal{G}_1)$
shows that $\Phi(S_i^0,S_j^1)\subset S_{i+j}^1$. 
From \cite{ver}, we also have
$$\left\{
\begin{aligned}
&\Phi(S_i^0,S_i^0)\subset S_{i+j+1}^0\text{ if }i+j<n\\
&\Phi(X_i,X_{n-i})=(-1)^i\alpha\ X_n\text{ where }\alpha=0\text{ if $n$ is even}\\
\end{aligned}
\right.$$
\end{proof}

This lead us to the study of the 2-cocycles of $L_{n,m}$.

\begin{prop} Let be $\Psi\in Z_0^2(L_{n,m},L_{n,m})$, such that $\mu_0+\rho_0+\Psi$
is a nilpotent Lie superalgebra,
then $\Psi$ admits the following decomposition $\Psi=\psi+\rho+b$ with
$$\begin{aligned}
\psi&\in Hom (\mathcal{G}_0\wedge \mathcal{G}_0,\mathcal{G}_0)\cap Z^2(L_{n},L_{n})\\
\rho&\in Hom(\mathcal{G}_0\otimes\mathcal{G}_1,\mathcal{G}_1)\cap Z^2(L_{n,m},L_{n,m})\\
b&\in Hom (\mathcal{G}_1\vee\mathcal{G}_1,\mathcal{G}_0)\cap Z^2(L_{n,m},L_{n,m})
\end{aligned}$$
\end{prop}

\begin{proof}
It is clear that $\Psi$ can be decomposed into a sum of three homogeneous maps :
$$\left\{
\begin{aligned}
\psi&\in Hom (\mathcal{G}_0\wedge \mathcal{G}_0,\mathcal{G}_0)\\
\rho&\in Hom(\mathcal{G}_0\otimes\mathcal{G}_1,\mathcal{G}_1)\\
b&\in Hom (\mathcal{G}_1\vee\mathcal{G}_1,\mathcal{G}_0)
\end{aligned}
\right.$$

As $\Psi\in Z_0^2(L_{n,m},L_{n,m})$, we have :
$$\left\{
\begin{aligned}
&\psi(\mu_0(g_i,g_j),g_k)-\psi(\mu_0(g_i,g_k),g_j)+\psi(\mu_0(g_j,g_k),g_i)-\\
&\quad \mu_0(g_i,\psi(g_j,g_k))+\mu_0(g_j,\psi(g_i,g_k))-\mu_0(g_k,\psi(g_i,g_j))=0,\\
&\\
&\rho(\mu_0(g_i,g_j),h_t)+\rho(g_j,\rho_0(g_i,h_t))-\rho(g_i,\rho_0(g_j,h_t))-\\
&\quad \rho_0(g_i,\rho(g_j,h_t))+\rho_0(g_j,\rho(g_i,h_t))-\rho_0(h_t,\psi(g_i,g_j))=0,\\
&\\
&b(\rho_0(g_i,h_t),h_r)+b(\rho_0(g_i,h_r),h_t)-\mu_0(g_i,b(h_t,h_r))=0,\\
&\\
&\rho_0(b(h_r,h_s),h_t)+\rho_0(b(h_t,h_s),h_r)+\rho_0(b(h_t,h_r),h_s)=0.
\end{aligned}
\right.$$
where $g_i$ is an even element and $h_j$ an odd element of $L_{n,m}$.
This prove that $\psi$ is a cocycle of the filiform Lie algebra $L_n$.
As $\mu_0+\rho_0+\Psi$ is nilpotent, $\psi(X_i,X_j)$ has no component on  $X_0$.
This implies that $\rho_0(h_t,\psi(g_i,g_j))=0$ and
$$\left\{\begin{aligned}
&\rho(\mu_0(g_i,g_j),h_t)+\rho(g_j,\rho_0(g_i,h_t))-\rho(g_i,\rho_0(g_j,h_t))-\\
&\quad \rho_0(g_i,\rho(g_j,h_t))+\rho_0(g_j,\rho(g_i,h_t))=0\\
&\\
&b(\rho_0(g_i,h_t),h_r)+b(\rho_0(g_i,h_r),h_t)-\mu_0(g_i,b(h_t,h_r))=0\\
&\\
&\rho_0(b(h_r,h_s),h_t)+\rho_0(b(h_t,h_s),h_r)+\rho_0(b(h_t,h_r),h_s)=0
\end{aligned}
\right.$$
This prove that both maps $\rho$ and $b$ are cocycles.
\end{proof}
We are reduced to study each space associated to the decomposition of $\Psi$.

\subsection{Cocycles of $Hom (\mathcal{G}_0\wedge \mathcal{G}_0,\mathcal{G}_0)$}
$\ $

Let $\psi$ be a 2-cocycle of $L_{n,m}$ belonging to  $Hom (\mathcal{G}_0\wedge \mathcal{G}_0,\mathcal{G}_0)$.
Then $\psi$ is a 2-cocycle of the Lie algebras $L_n$. From \cite{ver}, these cocycle are written as a linear sum of the following cocycles :

Let $(k,s)$ be a pair of integers such that $1\leq k\leq n-1$,
$2k\leq s\leq n$. There exists an  unique cocycle of $L_n$
satisfying  :
$$\Psi_{k,s}(X_i,X_{i+1})=\left\{
\begin{aligned} 
&X_s\text{ if }i=k\\
&0 \text{ otherwise }\end{aligned}
\right.$$
and
$$\Psi_{k,s}(X_0,X_i)=0\quad 1\leq i\leq n$$
is given for $i<j$ by 
$$\begin{aligned}
\Psi_{k,s}(X_i,X_j)&=(-1)^{k-i}C^{k-j}_{j-k-1} (ad\ X_0)^{i+j-2k-1} X_s\text{ if }k-i\leq j-k-1\\
\Psi_{k,s}(X_i,X_j)&=0\text{ otherwise }
\end{aligned}$$

The cocycles $\psi_{i,j}$ give the nilpotent deformations of the filiform Lie algebra $L_n$.

\subsection{Cocycles of $Hom (\mathcal{G}_0\oplus \mathcal{G}_1,\mathcal{G}_1)$}
$\ $
These cocycles are described in the following proposition :

\begin{prop}
For $1\leq k\leq n$ and $1\leq s\leq m$, there exists an unique cocycle $\rho_{k,s}$ of 
$Hom(\mathcal{G}_0\otimes\mathcal{G}_1,\mathcal{G}_1)\cap Z^2(L_{n,m},L_{n,m})$
such that :
$$\rho_{k,s}(X_i,Y_1)=\left\{
\begin{aligned}
&Y_s\text{ if }i=k\\
&0\text{ otherwise }
\end{aligned}
\right.$$
and
$$\rho_{k,s}(X_0,Y_i)=0\quad 1\leq i\leq m$$
It satisfies :
$$\rho_{k,s}(X_{j},Y_r)=
\left\{
\begin{aligned}
&(-1)^{k-j}C_{r-1}^{k-j}\ Y_{s+r-1-j+k}\text{ if } k-r+1\leq j\leq k\\
&0\text{ otherwise }
\end{aligned}
\right.$$
for $1\leq r\leq m$, where $C_{r-1}^{k-j}$ are the binomial coefficients.
\end{prop}

\begin{proof}

Let $\rho_{k,s}$ be a cocycle such that :
$$\rho_{k,s}(X_i,Y_1)=\left\{
\begin{aligned}
&Y_s\text{ if }i=k\\
&0\text{ otherwise }
\end{aligned}
\right.$$
and
$$\rho_{k,s}(X_0,Y_i)=0\quad 1\leq i\leq m$$
 $\rho_{k,s}$ must satisfy
$$\begin{aligned}
\rho_{k,s}([X_i,X_j],Y_r)&+\rho_{k,s}(X_j,[X_i,Y_r])-\rho_{k,s}(X_i,[X_j,Y_r])-\\
&[X_i,\rho_{k,s}(X_j,Y_r)]+[X_j,\rho_{k,s}(X_i,Y_r)]=0
\end{aligned}$$
then by induction on $r$ and $j$ we prove that
$$\rho_{k,s}(X_{j},Y_r)=
\left\{
\begin{aligned}
&(-1)^{k-j}C_{r-1}^{k-j}\ Y_{s+r-1-j+k}\text{ if } k-r+1\leq j\leq k\\
&0\text{ otherwise }
\end{aligned}
\right.$$
for $1\leq r\leq m$.
\end{proof}

\begin{prop}
Let $\{X_0,X_1,\dots,X_n,Y_1,Y_2,\dots ,Y_m\}$ be an adapted base of $L_{n,m}$.
The bilinear mapping $\varrho_{i,j}$ with $1\leq i,j\leq m$ defined by :
$$\left\{
\begin{aligned}
&\varrho_{i,j}(X_0,Y_i)=Y_j,\\
&\varrho_{i,j}(X_p,Y_k)=\varrho_{i,j}(X_p,Y_k)=\varrho_{i,j}(Y_p,Y_k)=0,\quad p\ne 0.
\end{aligned}
\right.$$ 
are cocycles.
\end{prop}
The prove is obvious.

\begin{thm}
The family of cocycles $\varrho_{i,j}$ and $\rho_{k,s}$
with $1\leq i\leq j\leq m$ and $1\leq k\leq n$, $1\leq s\leq m$
form a basis of  $Z_0^2(L_{n,m},L_{n,m})\cap Hom(\G_0\otimes\G_1,\G_1)$.
\end{thm}

\begin{proof}
Let $\rho$ be a cocycle of $Z_0^2(L_{n,m},L_{n,m})\cap Hom(\G_0\otimes\G_1,\G_1)$,
such that \\
$\rho(X_0,Y_j)=0$ for $1\leq j\leq m$. We can prove by induction on $j$ that
if $\rho(X_j,Y_1)=0$ for $1\leq j\leq n$ then $\rho\equiv 0$.

We can assume that $\rho(X_0,Y_j)=0$, if not, we consider the cocycle $\rho_1=\rho-\sum_{i,k}a_{i,k}\varrho_{i,k}$
such that $\rho_1(X_0,Y_j)=0$.
It is easy to see that there exists a linear combination of $\rho_{i,j}$ such that \\
$\rho'=\rho-\sum_{j=1}^m r_{i,j}\rho_{i,j}$ satisfies  $\rho'(X_i,Y_1)=0$. Using the previous
paragraph, we have that  $\rho'\equiv 0$ we deduce that $\rho=\sum_{j=1}^m r_{i,j}\rho_{i,j}$,
and if $\rho(X_0,Y_j)$ was not zero, $\rho$ will be 
$$\rho=\sum_{i=1}^n\sum_{j=1}^m t_{i,j}\ \rho_{i,j}+\sum_{k=1}^m \sum_{r=1}^{m} s_{k,r} \varrho_{k,r}$$
this prove that $\rho_{i,j}$ and $\varrho_{k,r}$ are generator. As they are linearly independent,
we have a base.
\end{proof}

\subsection{Cocycles of $Hom (\mathcal{G}_1\vee \mathcal{G}_1,\mathcal{G}_0)$}
$\ $

In this case, we will not give a basis for this cocycles,
but we will give the dimension of this space.

Let $b$ be a cocycle in $Hom (\mathcal{G}_1\vee\mathcal{G}_1,\mathcal{G}_0)\cap Z^2(L_{n,m},L_{n,m})$.
Then $b$ has to verify the two conditions :
\begin{equation}
\label{r1}
b(\rho_0(g_i,h_t),h_r)+b(\rho_0(g_i,h_r),h_t)-\mu_0(g_i,b(h_t,h_r))=0
\end{equation}
and
\begin{equation}
\label{r2}
\rho_0(b(h_r,h_s),h_t)+\rho_0(b(h_t,h_s),h_r)+\rho_0(b(h_t,h_r),h_s)=0
\end{equation}
where $g_i\in L_{n,m}^0$ and $h_t,h_r,h_s\in\ L_{n,m}^1$.

Now we will focus or work on relation (\ref{r1}). Note that if 
$g_i$ is linearly independent of $X_0$, then (\ref{r1}) is satisfied.
We suppose that $g_i=X_0$. Consider the adapted basis
$\{ X_0,X_1,\dots ,X_n,Y_1,\dots, Y_m\}$ of $L_{n,m}$ then (\ref{r1}) 
is written :
\begin{equation}
\label{r3}
\mu_0(X_0,b(Y_t,Y_r))=b(Y_{t+1},Y_r)+b(Y_{r+1},Y_t)
\end{equation}
for $1\leq r,t\leq m-1$.

\begin{lem}
\label{lbn}
Let $b$ be a symmetric bilinear mapping satisfying (\ref{r3}), such that
$$b(Y_i,Y_i)=0\text{ for }1\leq i\leq m;$$
then $b$ is null.
\end{lem}

\begin{proof}
Let us prove that $b(Y_i,Y_{i+k})=0$ for every $k$.
For $k=0$ we have :
$$b(Y_i,Y_{i})=0$$
Suppose that the relation is true up to $k$. For  $k+1$ we have :
$$\begin{aligned}
\ \mu_0(X_0,b(Y_i,Y_{i+k}))&=b(Y_{i+1},Y_{(i+1)+(k-1)})+b(Y_i,Y_{i+k+1})\\
0&=0+b(Y_i,Y_{i+k+1}).
\end{aligned}$$
then $b(Y_i,Y_{i+k})=0$ for all integer $k$ and $i$.
\end{proof}

This lemma shows that a symmetric bilinear map $b$ satisfying (\ref{r3}) can be
defined only by the value $b(Y_i,Y_i)$ with $1\leq i\leq m$.
Relation (\ref{r2}) implies that $\rho_0(Y_i,b(Y_i,Y_i))=0$, if
$b(Y_i,Y_i)=a_i\ X_0+\dots$, we have that $a_i.\rho(Y_i,X_0)=0$.
This implies that $a_i=0$, $1\leq i< m$.
Suppose that $a_m\ne 0$, then the
relation (\ref{r2}) implies that
$\rho_0(Y_1,b(Y_m,Y_m))=0$, then $a_m\ Y_2=0$ and $a_m=0$.
This prove that $b(Y_i,Y_i)$ does not have a component on $X_0$
for every $i$.

We define the vector space $E$ of the symmetric bilinear
maps satisfying the relation (\ref{r3}) such that $\forall b\in E, b(Y_m,Y_m)\in Vect(X_n)$. 

\begin{prop}
\label{fks}
The symmetric bilinear maps $f_{p,s}$   
with $1\leq s\leq n$ and $1\leq p\leq m-1$ defined 
for $1\leq i\leq p\leq j\leq m$  by :
$$f_{p,s} (Y_i,Y_j)=\frac{(-1)^{p-i}}{2}\left\lgroup C_{j-p}^{p-i}+C_{j-p-1}^{p-i-1}\right\rgroup X_{s-2p+i+j}$$
with convention $C_{-1}^{-1}=1$ and $0$ otherwise; and
$$f_{m,n}(Y_i,Y_j)=\left\{
\begin{aligned}
&X_s \text{ if }i=j=m\\
&0 \text{ otherwise}
\end{aligned}
\right.$$
form a basis of $E$.
\end{prop}

\begin{proof}
Let $b\in E$ be a symmetric bilinear map. 
It is easy to see that there exists coefficient
$a_{p,s}\in\mathbb{C}$ such that
\label{bsum}
$$b(Y_i,Y_i)-\sum_{p=1}^{m}\sum_{s=1}^n a_{p,s}\ f_{p,s}(Y_i,Y_i)=0$$
for $1\leq i\leq m$, where  $f_{m,i}=0$ for $1\leq i\leq n-1$.
Using lemma \ref{lbn} we deduce that this equality
vanishes for every pair $(Y_i,Y_j)$. This proves that
$$b=\sum_{p=1}^{m}\sum_{s=1}^n a_{p,s}\ f_{p,s}$$
Using the fact that $f_{p,s}(Y_p,Y_p)=X_s$ and $f_{p,s}(Y_i,Y_i)=0$ if
$i\ne p$, the family $\{f_{p,s}\}$ is free.
\end{proof}

\begin{prop}
\label{structcf}
The space
$Z_0^2(L_{n,m},L_{n,m})\cap Hom(\G_1\vee \G_1,\G_0)$
is the subspace of $E$ defined by
$$(5)\left
\{ 
\begin{aligned}
&f\in E,\\
&\mu_0(X_0,f(Y_i,Y_m))=f(Y_{i+1},Y_m)\\
&\text{ for }1\leq i\leq m-1
\end{aligned}
\right \}$$

\end{prop}

\begin{proof}
A cocycle $f$ satisfies the two 
relations (\ref{r1}) and (\ref{r2}). A consequence of this
is that $f(Y_m,Y_m)=a_{m,n}\ X_n$.  
We deduce that $f\in E$.
To satisfy relation (\ref{r2}), $f$ has to satisfy the relation
$$\mu_0(X_0,f(Y_i,Y_m))=f(Y_{i+1},Y_m)$$
for $1\leq i\leq m-1$.

Such a map $f$ does not have a component on $X_0$ in its image, hence
$[Y_i,b(Y_j,Y_k)]=0$ for $1\leq i,j,k\leq m$. This prove that relation (\ref{r1}) is satisfied
and that every map $f$ satisfying (\ref{r2}) is a cocycle.
\end{proof}

Consequence : 
$$\dim Z_0^2(L_{n,m},L_{n,m})\cap Hom(\G_1\vee \G_1,\G_0)\leq \dim E=n.m-n+1$$ 

The maps $f_{p,s}$ are not always cocycles. Let  $b_{p,s}^{(\alpha )}\in E$ be
 $$b_{p,s}^{(\alpha )}=f_{p,s}+\sum_{k=1}^{m-p} \alpha_{p,s}^k f_{p+k,s+2k}$$
where $\alpha_{p,s}^k\in\mathbb{C}$. Let $A_{p,s}$ be
the set of sequences $(\alpha )=(\alpha_{p,s}^1,\alpha_{p,s}^2,\dots \alpha_{p,s}^{m-p})$ such
that $b_{p,s}^{(\alpha)}$ is a cocycle.
Then  $(\alpha_{p,s}^1,\alpha_{p,s}^2,\dots \alpha_{p,s}^{m-p})$ is a solution
of equation (5). Remark that for some pair $(p,s)$, $A_{p,s}$ can be empty.

\begin{lem}
The family of cocycles $b_{p,s}^{(\alpha)}$ with $(\alpha)\in\cup_{p,s} A_{p,s}$
spans 
$$Hom (\mathcal{G}_1\vee\mathcal{G}_1,\mathcal{G}_0)\cap Z^2(L_{n,m},L_{n,m})$$
\end{lem}

\begin{proof}

Using theorem \ref{structcf} every cocycle of  $Hom(\G_1\vee \G_1,\G_0)$
is given by :
$$f=\sum_{p=p_0}^{m}\sum_{s=s_0}^n a_{p,s}\ f_{p,s}$$
with $a_{p_0,s_0}\ne 0$.

We can write  $f$ like :
$$f=a_{s_0,p_0}\ b_{s_0,p_0}^{(\alpha )}+R$$
where $b_{s_0,p_0}(Y_i,Y_j)$ has a component on $X_{s_0-2p_0+i+j}$ and $R(Y_i,Y_j)$ does not have
any component on $X_{s_0-2p_0+i+j}$.

As $f$ is a cocycle, we have :
$$\begin{aligned}
\ [X_0,f(Y_i,Y_m)]&=f(Y_{i+1},Y_m)\\
a_{s_0,p_0}\ [X_0,b_{s_0,p_0}(Y_i,Y_m)]+[X_0,R(Y_i,Y_j)]&=a_{s_0,p_0}\ b_{s_0,p_0}(Y_{i+1},Y_m)+R(Y_{i+1},Y_m)\\
\end{aligned}$$

Let us consider the component on $X_{s_0-2p_0+i+1+m}$, we have :
$$\begin{aligned}
a_{s_0,p_0}\ [X_0,b_{s_0,p_0}(Y_i,Y_m)]&=a_{s_0,p_0}\ b_{s_0,p_0}(Y_{i+1},Y_m)\\
[X_0,b_{s_0,p_0}(Y_i,Y_m)]&=b_{s_0,p_0}(Y_{i+1},Y_m)\text{ because }a_{s_0,p_0}\ne 0\\
\end{aligned}$$
This proves that $b_{s_0,p_0}$ is a cocycle, as  $f-a_{s_0,p_0}\ b_{s_0,p_0}$.
Using the cocycle $f-a_{s_0,p_0}\ b_{s_0,p_0}$,
we prove by induction that $f$ is given by a linear combination of the cocycles  $b_{p,s}^{(\alpha)}$.

\end{proof}

The cocycles $b_{p,s}^{(\alpha)}$ are not linearly independent. Therefore, for a non empty set
$A_{p,s}$, we consider the smallest cocycle $b_{p,s}^0$ in the sense that $b_{p,s}^0$ cannot
be written 
$$b_{p,s}^0=b_{p,s}^{(\alpha_1)}+a\ b_{k,s}^{(\alpha_2)}$$
with $b_{k,s}^{(\alpha_2)}$ non zero, $a\in \mathbb{C}^*$ and $k>p$.

\begin{lem}
\label{eci}
For any non empty set $A_{p,s}$ there exists an unique cocycle $b_{p,s}^0$.
\end{lem}

\begin{proof}
Let be $A_{p,s}\ne \emptyset$ and $b_{p,s}^{(\alpha)}$ be a non zero cocycle.
If we can decompose $b_{p,s}^{(\alpha)}$, we have the smallest cocycle, if not we choose 
the smallest integer $k_0$, $p<k_0+p\leq m$ such that
$$b_{p,s}^{(\alpha)}=b_{p,s}^{(\alpha_0)}+\gamma_{p,s}^{k_0}\ b_{p+k_0,s+2k_0}^{(\alpha_k)}$$
If $b_{p,s}^{(\alpha_0)}$ is indecomposable, we stop.
If not, we have
$$b_{p,s}^{(\alpha_0)}=b_{p,s}^{(\alpha_1)}+\gamma_{p,s}^{k_1}\ b_{p+k_1,s+2k_1}^{(\alpha_2)}$$
with $k_1>k_0$, this sequence is increasing and has an upper bound, therefore it exists $k_r$ such that
$b_{p,s}^{(\alpha_r)}$ is indecomposable.

To proof the uniqueness, suppose that
$b_{p,s}^{(\alpha_1)}$ and $b_{p,s}^{(\alpha_2)}$ are smallest. Then 
$b_{p,s}^1-b_{p,s}^2=b$ is a cocycle. We have 
$b_{p,s}^{(\alpha_1)}=b_{p,s}^{(\alpha_2)}+a\ b_{k_0,r}^{(\alpha_3)}$ with $k_0$ the
smallest integer $k$ such that $f_{k,r}$ is in  $b$.
 As $b_{p,s}^{(\alpha_1)}$ is the smallest, we must have
$a\ b_{k_0,r}=0$ and then  $b_{p,s}^{(\alpha_1)}=b_{p,s}^{(\alpha_2)}$.
This proves the uniqueness.
\end{proof}

\begin{thm}
\label{thm1}
The cocycles' family  $b_{p,s}^0$ with $(p,s)$ such that $A_{p,s}\neq\emptyset$ is a basis of 
$$Hom (\mathcal{G}_1\vee\mathcal{G}_1,\mathcal{G}_0)\cap Z^2(L_{n,m},L_{n,m})$$
\end{thm}

\begin{proof}
Let $(p,s)$ be such that $A_{p,s}\neq\emptyset$.
The cocycles $b_{p,s}^0$ span $Hom (\mathcal{G}_1\vee\mathcal{G}_1,\mathcal{G}_0)\cap Z^2(L_{n,m},L_{n,m})$
because every cocycle $b_{p,s}$ can be written as a linear sum of $b_{p+k,s+2k}^0$, $k\geq 0$.

Let be $a_{p,s}\in\mathbb{C}$ such that :
$$\sum_{p=1}^m \sum_{s=1}^n a_{p,s}\ b_{p,s}^0\equiv 0$$
Note that :
$$b_{p,s}^0(Y_i,Y_i)=\left\{
\begin{aligned}
&0\text{ if }i<p\\
&X_s\text{ if }i=p\\
&\alpha_{p,s}^k X_s\text{ if }i=p+k
\end{aligned}
\right.$$
We have for $Y_1$ :
$$\sum_{p=1}^m \sum_{s=1}^n a_{p,s}\ b_{p,s}^0(Y_{1},Y_{1})=\sum_{s=1}^n a_{1,s}\ X_s$$
then $a_{1,s}=0$ for $1\leq s\leq n$.
By induction on $p=1,2,\dots m$, every coefficient vanishes, and the $b_{p,s}^0$
are linearly independent.

\end{proof}

The theorem shows that the determination of a basis of 
$$Hom (\mathcal{G}_1\vee\mathcal{G}_1,\mathcal{G}_0)\cap Z^2(L_{n,m},L_{n,m})$$
is reduced to the case $A_{p,s}\neq \emptyset$.

\begin{prop}
\label{prop1}
The only pairs $(p,s)$ such that $f_{p,s}$ is a cocycle,
that is $(0,0,\dots )\in A_{p,s}$, are

\begin{itemize}
\item if $m$ is odd : 
$f_{\frac{m-1}{2},n}$ and $f_{p,s}$ with $2p=m-k$, $n-k-1\leq s\leq n$ for
$1\leq k\leq m-2$ and $k$ odd.

\item if $m$ is even :
$f_{p,s}$ with $2p=m-k$, $n-k-1\leq s\leq n$ for
$0\leq k\leq m-2$ and $k$ even.
\end{itemize}

\end{prop}

\begin{proof}

Let $f_{p,s}$ be a cocycle from the proposition.
If $f_{p,s}\in Z^2(L_{n,m},L_{n,m})$ then 
$$[X_0,f_{p,s}(Y_i,Y_m)]=f_{p,s}(Y_{i+1},Y_m)$$

Let $(p,s)$ be such that  $2p=m-k$ and $1\leq s\leq n-k-2$, we will proof that
$f_{p,s}$ is not a cocycle. We have
$$f_{p,s}(Y_1,Y_m)=\frac{(-1)^{p-1}}{2} (C_{m-p}^{p-1}+C_{m-p-1}^{p-2})X_{s+k+1}$$
$$f_{p,s}(Y_2,Y_m)=\frac{(-1)^{p-2}}{2} (C_{m-p}^{p-2}+C_{m-p-1}^{p-3})X_{s+k+2}$$
If $f_{p,s}$ is a cocycle, we have
$$(C_{m-p}^{p-1}+C_{m-p-1}^{p-2})=-(C_{m-p}^{p-2}+C_{m-p-1}^{p-3})$$
This implies
$$\begin{aligned}
C_{m-p}^{p-1}+&C_{m-p-1}^{p-2}=0\\
C_{m-p}^{p-2}+&C_{m-p-1}^{p-3}=0
\end{aligned}$$
Thus $C_{m-p}^{p-1}=0$, and
$2p>m+1$. 
We have $f_{p,s}(Y_1,Y_m)=0$ and $f_{p,s}(Y_{2p-m},Y_m)=(-1)^{m-1}\ X_s\neq 0$.
 This proofs that $f_{p,s}$ cannot be a cocycle.

\end{proof}

\begin{prop}
\label{prop2}
Let $q$ be such that  $1\leq q\leq\min\{ m-1,n-2\}$.
If $p$ satisfies $2+m+q-n\leq 2p\leq m-q+1$ then
$A_{p,s}$, with $s=n-m-q-1+2p$, is not empty.
\end{prop}

\begin{proof}
Let $q$ be such that  $1\leq q\leq\min\{ m-1,n-2\}$, 
$s=n-m-q-1+2p$
and $p$ such that $2+m+q-n\leq 2p\leq m-q+1$.

Let's proof that there exists $(\alpha_{p,s}^1,\dots ,\alpha_{p,s}^q)\in A_{p,s}$ 
such that 
$$b_{p,s}=f_{p,s}+\sum_{k=1}^{q} \alpha_{p,s}^k\ f_{p+k,s+2k}$$
is a cocycle.

We have $b_{p,s}(Y_1,Y_m)\in \mathbb{C}.X_{n-q}$, then for $b_{p,s}$ to be a cocycle,
it must satisfy $q$ equations given by (5). Suppose $1\leq i\leq q$, if 
$$f_{p+k,s+2k}(Y_i,Y_m)=\frac{(-1)^{p+k-i}}{2} (C_{m-p-k}^{p+k-i}+C_{m-p-k-1}^{p+k-i-1})X_{n-q+i-1}$$
is vanishing then $C_{m-p-k}^{p+k-i}+C_{m-p-k-1}^{p+k-i-1}=0$ as $n-q+i-1\leq n$. This is
possible only if $p+k-i-1\leq m-p-k-1$. As $2+m+q-n\leq 2p\leq m-q+1$, there exists
a value of $p$ and $k$ such that $C_{m-p-k}^{p+k-i}+C_{m-p-k-1}^{p+k-i-1}\ne 0$.
 This proofs that in the $q$ linear equations given by
(5), every coefficient $\alpha_{p,s}^1,\dots ,\alpha_{p,s}^q$ appear.
This proof that  this system admits a solution.
\end{proof}

\begin{cor}
Suppose that $m\geq n$, $m=2t$. Then if
\begin{itemize}
\item $n=4s$ : $\dim Z_0^2(L_{n,m},L_{n,m})\cap Hom(\G_1\vee \G_1,\G_0)\geq t.n-2s^2+s$
\item $n=4s+1$ : $\dim Z_0^2(L_{n,m},L_{n,m})\cap Hom(\G_1\vee \G_1,\G_0)\geq t.n-2s^2$
\item $n=4s+2$ : $\dim Z_0^2(L_{n,m},L_{n,m})\cap Hom(\G_1\vee \G_1,\G_0)\geq t.n-2s^2-s$
\item $n=4s+3$ : $\dim Z_0^2(L_{n,m},L_{n,m})\cap Hom(\G_1\vee \G_1,\G_0)\geq t.n-2s^2-s$
\end{itemize}
For $m\geq n$, $m=2t+1$ then if
\begin{itemize}
\item $n=4s$ : $\dim Z_0^2(L_{n,m},L_{n,m})\cap Hom(\G_1\vee \G_1,\G_0)\geq (t+1).n-2s^2-s$
\item $n=4s+1$ : $\dim Z_0^2(L_{n,m},L_{n,m})\cap Hom(\G_1\vee \G_1,\G_0)\geq t.n-2s^2+2s+1$
\item $n=4s+2$ : $\dim Z_0^2(L_{n,m},L_{n,m})\cap Hom(\G_1\vee \G_1,\G_0)\geq (t+1).n-2s^2-4s-1$
\item $n=4s+3$ : $\dim Z_0^2(L_{n,m},L_{n,m})\cap Hom(\G_1\vee \G_1,\G_0)\geq t.n-2s^2+2$
\end{itemize}
\end{cor}

\begin{proof}
Using propositions \ref{prop1} and \ref{prop2}, we can compute
a lower bound of the number of non empty sets $A_{p,s}$. For each of this sets,
there exists a unique cocycle $b_{p,s}^0$ (see lemma \ref{eci})
which is a vector of the base of  
$$Hom (\mathcal{G}_1\vee\mathcal{G}_1,\mathcal{G}_0)\cap Z^2(L_{n,m},L_{n,m})$$
this is established in the theorem \ref{thm1}.
\end{proof}

I should like to thank Y. Khakimdjanov and M. Goze for support and numerous discussions.
\end{document}